\documentclass[a4paper,12pt]{amsart}
\usepackage{amsfonts}
\usepackage{amssymb}
\usepackage{amsfonts}
\usepackage{amssymb}
\usepackage{amsthm}
\usepackage{amsmath}
\usepackage{a4wide}

\numberwithin{equation}{section}

\theoremstyle{plain}
\newtheorem{teo}{Theorem}[section]
\newtheorem{lemma}[teo]{Lemma}
\newtheorem{prop}[teo]{Proposition}

\theoremstyle{definition}
\newtheorem{dfnz}[teo]{Definition}

\theoremstyle{remark}
\newtheorem{rem}[teo]{Remark}

\newtheorem{es}[teo]{Example}

\def\R{\mathbb R}

\def\NN{\mathbb N}
\def\composed{\circ}
\def\comp{\composed}
\def\eps{\varepsilon}

\def\dert{\partial_t}
\def\Derpar#1 { \frac{\partial~ } {\partial {#1} }} 
\def \derpar #1 #2 { \frac {\partial{#2}}{\partial {#1}}}
\def\ders{\partial_s}
\def\pol{{\mathfrak{p}}}
\def\qol{{\mathfrak{q}}}
\def\fol{{\mathfrak{f}}}

\def\composed{\circ}

\def\DG{{\mathcal G}_k}
\def\HHH{{\mathrm H}}

\def\bbigstar{\operatornamewithlimits{\text{\Large{$\circledast$}}}}
\def\EEE{{\mathrm E}}
\def\BBB{{\mathrm B}}

\def\ZZZ{{\mathrm Z}}

\newcommand{\Eul}         {{{\EEE}^\varepsilon}}
\renewcommand{\k}         {\kappa}

\renewcommand{\L}       {\mathrm{L}}

\begin{document}

\title[Singular Perturbations of Mean Curvature Flow]
{Singular Perturbations of Mean Curvature Flow}

\author[Giovanni Bellettini]{Giovanni Bellettini}
\address[Giovanni Bellettini]{Dipartimento di Matematica, Univ. Roma ``Tor 
Vergata'',
  Roma, Italy}
\email[G.~Bellettini]{belletti@mat.uniroma2.it}

\author[Carlo Mantegazza]{Carlo Mantegazza}
\address[Carlo Mantegazza]{Scuola
  Normale Superiore, Pisa, 56126, Italy}
\email[C.~Mantegazza]{mantegazza@sns.it}

\author[Matteo Novaga]{Matteo Novaga}
\address[Matteo Novaga]{Dipartimento di Matematica, Univ. Pisa, Pisa, 
56127, Italy}
\email[M.~Novaga]{novaga@dm.unipi.it}

\date{\today}

\begin{abstract}
We introduce a regularization method for mean
curvature flow of a submanifold of arbitrary codimension in the
Euclidean space, through higher order equations.\\ 
We prove that the regularized problems converge to the mean curvature flow
for all times before the first singularity.
\end{abstract}

\maketitle

\tableofcontents

\section{Introduction}
It is well known 
that a smooth compact submanifold of the Euclidean space, flowing by
mean curvature, develops singularities in finite time.
This is a common aspect of geometric evolutions, and 
motivates the study of the flow
past singularities. Concerning the mean curvature motion, 
several notions of weak solutions have been proposed, after the
pioneering work of Brakke~\cite{brakke}, see for
instance~\cite{altawa,ambson,ambson97,BeNo:97,BeOrSm:03,cgg,
degio4,es,ilman2,Il:93g,ilman5,ilman1,soner1}.
We recall that some of these solutions may differ, in 
particular in presence of the so--called fattening phenomenon
(see for instance~\cite{BaSoSo:93}).\\ 
Following a suggestion of De Giorgi in~\cite{degio5}, 
we introduce and study a regularization of mean curvature 
flow with a singular perturbation of higher order, which could lead to a
new definition of generalized solution in any dimension and
codimension.

Let us state our main result.\\
Let $\varphi:M\to\R^{n+m}$ be a smooth compact $n$--dimensional 
immersion in $\R^{n+m}$. 
For $k>[n/2]+2$ (where $[n/2]$ denotes the integer part of $n/2$) 
and $\varepsilon>0$ we consider the functional
$$
\DG^\varepsilon(\varphi) = \int_M
\left( 1+\varepsilon\vert A^k\vert^2\right)\,d\mu\,,
$$
where $\mu$ is the canonical volume measure associate with the metric
$g$ induced on $M$ via the immersion $\varphi$. With $A^k$ we denote
the $k$--differential in $\R^{n+m}$ of the function $A^M$ given by
$$
A^M(x)=\frac{\vert x\vert^2-[d^M(x)]^2}{2}\,,
$$
where $[d^M]^2$ is the square of the
distance function from $\varphi(M)$, which is smooth in a
neighborhood of a point of the submanifold without
self--intersections. Since locally on 
$M$ every immersion is an embedding, we can define $A^k$ also at such
points. More precisely, the tensor $A^k$ is defined as
$$
A^k_{i_1\dots i_k}=\frac{\partial^k A^M}{\partial x_{i_1}\dots\partial
  x_{i_k}}
$$
for every $k$--uple of indexes $i_1,\dots, i_k\in\{1,\dots,n+m\}$.\\
We remark (see~\cite[Prop.~2.2 and Cor.~2.4]{mantemin1}) that the
tensors $A^k$ and $\nabla^{k-3}\BBB$, where $\BBB$ is the second
fundamental form of $M$ and $\nabla$ is the covariant derivative
associated with the induced metric $g$, are strictly related, hence, in 
a way the functional $\DG^\varepsilon$ is a perturbation of the area
functional by a term containing the squares of the high order 
derivatives of the curvatures of $M$.

c
By means of Theorem~4.5 and Theorem~5.9 in~\cite{ambman1} and the
results of~\cite{mantemin1}, the {\em gradient flow} associated 
with the functional 
$\DG^\varepsilon$ is given by the PDE system 
\begin{equation}\label{tttt1}
\frac{\partial\varphi^\varepsilon}{\partial t}=
\HHH+ 2\varepsilon k(-1)^{k}
\Bigl(\overset{\text{$( k -2)$--times}}{\overbrace{\Delta^{M}\comp
\Delta^{M}\comp\,\dots\,\comp\Delta^{M}}}\,\,\HHH\Bigr)^\perp
+\varepsilon\,{\mathrm {LOT}}
\end{equation}
where $\HHH$ is the mean curvature vector and ${\mathrm {LOT}}$ denotes  
terms of lower order in the curvature and its derivatives.\\
We can see then that~\eqref{tttt1} is a {\em singular perturbation} of 
the 
mean
curvature flow, and coincides with it when $\eps=0$.
In~\cite{mantemin1} (see also~\cite{emintes} and~\cite{mant5}) it is
proved that for every $\eps>0$ the system in~\eqref{tttt1} 
admits a unique smooth solution 
defined 
{\em for all times};
we are then interested in the convergence to the mean curvature flow 
when $\eps\to0$.\\
Our main result is the following.
\begin{teo}\label{mainth} 
Let $\varphi_0:M\to\R^{n+m}$ be a smooth immersion of a compact
$n$--dimensional manifold without boundary.
Let $T_{\mathrm {sing}}>0$ be the first singularity time of 
the mean curvature flow $\varphi:M\times[0,T_{\mathrm {sing}})\to\R^{n+m}$ 
of $M$. For any $\eps>0$ let
$\varphi^\eps:M\times[0,+\infty)\to\R^{n+m}$ 
be the flows associated with the functionals $\DG^\eps$, with
$k>[n/2]+2$, all starting from the same initial immersion $\varphi_0$.
Then the maps $\varphi^\eps$ converge locally in
$C^\infty(M\times[0,T_{\mathrm {sing}}))$ to the map $\varphi$, as $\eps\to0$.
\end{teo}

\begin{es}
In case of immersed plane curves $\gamma: S^1 \to \R^2$ ($n=m=1$)
the simplest choice is $k=3$. 
Since it turns out that $\vert A^3\vert^2 = 3 \kappa^2$, where
$\kappa$ is the curvature of $\gamma$, in this simple case
the approximating functionals 
read as
\begin{equation*}
\int_{\gamma} \left(1 + \eps
\kappa^2\right)\, ds
\end{equation*}
where $s$ is the arclength parameter, and we have
replaced $3\varepsilon$ with $\varepsilon$. The 
regularized system 
which approximate the curve shortening flow is then
\begin{equation}
\frac{\partial \gamma}{\partial t} =
\left(\kappa - 2 \eps \partial_s^2 \kappa - \eps \kappa^3 \right) \nu,
\end{equation}
where $\nu$ is a suitable choice of the normal unit vector to the curve.
\end{es}
The crucial point in order to prove Theorem~\ref{mainth} is to 
obtain $\varepsilon$--independent estimates of the curvature and its
derivatives in order to gain sufficient compactness properties. We get
these by computing the evolution equations satisfied by the
$L^2$ norms of the derivatives of the second fundamental form of the
flowing manifolds, and by estimating via Gagliardo--Nirenberg
interpolation inequalities.\\
At the present moment we are not able to 
characterize 
the limit of the approximating flows 
after the first singularity, as the proof of
Theorem~\ref{mainth} relies heavily on the smoothness of the mean
curvature flow in the time interval of existence. Our goal would be
to provide some limit flow defined {\em for all times}, thus providing a 
new weak definition of solution in any dimension and codimension.\\
We mention the simplest open problem in defining a
limit flow after the first singularity.\\
It is well known (Gage--Hamilton~\cite{gage,gaha1} and
Huisken~\cite{huisk1}) that a convex curve in the plane (or 
hypersurface in
$\R^{n+1}$) moving by mean curvature shrinks to a point in finite time,
becoming exponentially round. In this 
case we expect that the approximating flows 
converge (in a way to be made precise) to such point for
every time after the extinction one.

The plan of the paper is the following. In Section~\ref{secnot} we give
some notation and we recall the relations between the squared
distance function and the second fundamental 
form and its covariant derivatives. In Section~\ref{seccurves}, in
order to make the line of proof clearer, we work out in detail the
$\eps$--independent estimates in the simplest case of plane
immersed curves; also in this special case, the result appears
to be nontrivial. In Sections~\ref{secgen},~\ref{evolgeosec}
and~\ref{secepsindep} we consider the general case of a
$n$--dimensional submanifold of $\R^{n+m}$.
Section~\ref{seclim} is devoted to show Theorem~\ref{mainth}.

We remark here (but we will not discuss such an extension in this paper) 
that 
our method works in general for any geometric evolution of
submanifolds in a Riemannian manifold till the first singularity time,
even when the equations are of high order (like, for instance, in the
Willmore flow, see~\cite{kuschat1,kuschat2,willmore}), 
choosing a regularizing term of appropriately high order.

Finally, it should be noted that, looking at the evolution 
equation~\eqref{tttt1}, 
these perturbations of the mean curvature flow could be 
considered, in the framework of geometric evolution problems,  as an 
analogue of the so--called {\em vanishing viscosity method}.
Indeed, we perturb the mean curvature flow equation with a regularizing
higher order
term multiplied by a small parameter $\varepsilon>0$. The lower
order terms, denoted by ${\mathrm{LOT}}$, which appear in~\eqref{tttt1}
are due to the fact that we actually perturb the area functional and
not directly the evolution equation.
However, the analogy with the classical viscosity 
method cannot be pushed too far. For instance, because of the
condition $k>[n/2]+2$, our regularized
equations are of order not less than four (precisely at least four for
evolving curves, at least  six for evolving 
surfaces). Moreover, as the Laplacians appearing in 
equation~\eqref{tttt1} are relative to the induced metric, the system
is quasilinear  and the lower order terms are non linear
(polynomial).

\section{Notation and Preliminaries}\label{secnot}

We denote with $e_1,\dots,e_{n+m}$ the canonical basis of
$\R^{n+m}$ and with $\langle\,\,,\,\,\rangle$ its standard scalar product.\\
We let $M\subset\R^{n+m}$ be a smooth, compact,
$n$--dimensional, regular submanifold without boundary, then 
$T_xM$, $N_xM\subset\R^{n+m}$ are, respectively, the tangent space
and the normal space to $M$ at $x\in M\subset\R^{n+m}$.

The {\em distance function} $d^M$ and the
{\em squared distance function} 
$\eta^M$ from $M$ are given by
$$
d^M(x)=\inf_{y\in M}\vert x-y\vert\qquad\text{ and }\qquad 
\eta^M(x)=[d^M(x)]^2
$$
for any $x\in\R^{n+m}$ (we will drop the superscript $M$ when no ambiguity
is possible). In this section we recall some facts from~\cite{ambman1}
and~\cite{mantemin1} about the distance function and the relations
between the high derivatives of $\eta^M$ and the second fundamental
form of $M$.

When $M$ is embedded, there exists an open
neighborhood $\Omega\subset\R^{n+m}$ of $M$ such that 
$d^M$ is smooth in $\Omega\setminus M$ and $\eta^M$ is smooth in
the whole of $\Omega$. If $M$ is only immersed, at any point of $M$ we 
consider the distance (we still use the symbols $d^M$ and $\eta^M$ for 
simplicity) from 
an embedded image of a suitable neighborhood of the point; in this case 
the regularity properties of $d^M$ and $\eta^M$ hold in a neighborhood 
(still denoted by $\Omega$) of such an embedded image.
\\
Clearly, $\eta^M$ and $\nabla\eta^M(x)=0$ at every $x\in M$;
moreover, for every $x\in\Omega$ we have that $x-\nabla\eta^M(x)/2$ is
{\em the unique} point in $M$ of minimum distance from $x$ (the {\em
  projection} of $x$ on $M$), that we denote with $\pi^M(x)$.\\
Another nice property of the squared distance is that, for every 
$x\in M$ the Hessian matrix $\nabla^2\eta^M(x)$ is twice the matrix of
orthogonal projection onto the normal space $N_xM$. We will denote
respectively with $X^M$ and $X^\perp$ the orthogonal projections of a vector $X$
on the tangent and normal space of $M$.

Let $x\in M$ and $X,\,Y\in T_x M$, the vector valued 
second fundamental form of $M$ at the point $x$ is given by 
$$
\BBB(X,Y)=\Bigl(\frac{\partial Y}{\partial X}\Bigr)^\perp, 
$$
where we extended locally the two vectors $X$, $Y$ to 
tangent vector fields on $M$ (the derivative is
well defined since $X$ is a tangent vector at $x$).\\
If $\{\nu_\alpha\}_{\alpha=1,\dots, m}$ is a local basis of the normal
  bundle we have 
$$
\BBB(X,Y)=-\Bigl\langle\frac{\partial
  \nu_\alpha(x)}{\partial X},Y\Bigr\rangle \nu_\alpha(x)\,,
$$
where here and throughout all the paper we use the convention of summing
on repeated indices.\\
We will see $\BBB$ as a bilinear map from $T_xM\times T_xM$ to
$\R^{n+m}$, hence, as a family of $n+m$ 
bilinear forms $\BBB^k=\langle\BBB,e_k\rangle: T_xM\times
T_xM\to\R$. Moreover, we consider $\BBB$ acting also on vectors of
$\R^{n+m}$, not necessarily tangent, by setting 
$\BBB(V,W)=\BBB(V^M,W^M)\in N_xM\subset \R^{n+m}$ for every pair 
$V,\, W\in\R^{n+m}$. With such definition,
$\BBB^k_{ij}=\langle\BBB(e_i,e_j),e_k\rangle$.\\
It is well known that ${\BBB}$ is a symmetric bilinear form and its
trace is the {\em mean curvature} of components
${\HHH^k}=\BBB^k_{jj}$.

We introduce now the function
$$
A^M(x)=\frac{\vert x\vert^2-[d^M(x)]^2}{2}\,,
$$
smooth as $\eta^M$ in the neighborhood $\Omega$ of $M$, and we set
$$
A^M_{i_1 \dots i_k}(x)=\frac{\partial^kA^M(x)}{\partial x_{i_1}\dots\partial
  x_{i_k}}
$$
for the derivatives of $A^M$ at every point $x\in\Omega$.\\
The  following Proposition (see~\cite{ambman1} for the proof) shows the first 
connection between the second fundamental form and the function $A^M$
(or equivalently, the squared distance function).

\begin{prop} The following relations hold.
\begin{itemize}
\item For any $x\in\Omega$, the point $\nabla A^M(x)$ is the projection
point $\pi^M(x)$.
\item If $x\in M$, then $\nabla^2A^M(x)$ is the matrix of orthogonal projection on $T_xM$.
\item For every $x\in M$, 
\begin{align*}
\BBB^k_{ij}=&\,A^M_{ijs}(\delta_{ks}- A^M_{ks})\,, \\
A^M_{ijk}=&\,\BBB^k_{ij}+\BBB^i_{jk}+\BBB^j_{ki}\,,\\
\HHH^k=&\, A^M_{iik}\,.
\end{align*}
\end{itemize}
\end{prop}

We define now the $k$--derivative tensor $A^{k}(x)$ acting on 
$k$--uples
of vectors $X_i\in\R^{n+m}$, where $X_i=X_i^je_j$, as follows
$$
A^{k}(x)(X_1, \dots, X_k)=A^M_{i_1 \dots i_k}(x)X_1^{i_1}\dots
X_k^{i_k}\,,
$$
notice that the tensors $A^k$ are symmetric.\\
For notational simplicity, we drop the superscript $M$ on $A^k$;
for the
same reason, we also avoid to indicate the point $x\in M$ in the
sequel.\\
The tensors $A^k$ and $\nabla^{k-3}\BBB$ are strictly related by a
recurrence formula proved in~\cite[Prop.~2.2 and Cor.~2.4]{mantemin1}.

\begin{rem}
We underline here an important convention used in the paper.
Due to the high codimension, we will work with 
several tensors (like normal vector fields or the second
fundamental form $B$) taking values in $\R^{n+m}$; these
tensors will be considered as {\em families} of ${n+m}$ tensors with 
values in $\R$. With this convention, for instance, $\nabla\BBB$ means 
that we are considering the family of covariant derivatives of the
tensors $\BBB^1,\dots,\BBB^{n+m}$, one component at 
time. Unless otherwise specified, this convention
 will be used even also for {\em tangent} vector
fields, that is, when $X$ is a tangent vector field, $\nabla X$ is not the
covariant derivative of $X$ but the derivative of its components in
the basis of $\R^{n+m}$.
\end{rem}

In all the paper we write $T * S$, following
Hamilton~\cite{hamilton1}, to denote a tensor formed by contraction on 
some indices of the tensors $T$ and $S$ using the coefficients 
$g^{ij}$.\\
If $T_1, \dots , T_l$ are tensors
(here $l$ is not an index of the tensor $T$), 
with the symbol
$$
\bbigstar_{i=1}^l T_i
$$
we mean $T_1 * T_2 * \dots * T_l$\,.

\begin{dfnz}\label{polqol}
We use the symbol $\pol^s(\nabla^l\BBB)$ for a polynomial (with
the $*$ product) tensor with constant coefficients in the coordinate basis
$\partial\varphi/\partial{x_i}$, the second
fundamental form $\BBB$ and its derivatives up to the order $l$ at
most, such that every of its monomial is of the form
$$
\bbigstar_{k=1}^N\nabla^{j_k}\BBB
\qquad\text{ with $0 \leq j_k\leq l$ and $N\geq1$}
$$
or
$$
\bbigstar_{k=1}^N\nabla^{j_k}\BBB \frac{\partial\varphi}{\partial{x_i}}
\qquad\text{ with $0 \leq j_k\leq l$ and $N\geq1$}
$$
where, in both cases, the {\em rescaling order} $s$ equals
$$
s=\sum_{k=1}^N (j_k+1)\,.
$$
We use instead the symbol $\qol^s(\nabla^l\BBB)$ for a polynomial
of the kind $\pol^s(\nabla^l\BBB)$ such that 
the contraction with the metric is total, both in the covariant and
in the $\R^{n+m}$--indices.\\
As the contraction in the ambient space $\R^{n+m}$ ``cancels'' all the
basis elements $\partial\varphi/\partial x_i$ appearing in
the formulae, it 
follows that 
every monomial of $\qol^s(\nabla^l\BBB)$ has the form
$$
\bbigstar_{k=1}^N\nabla^{j_k}\BBB
\qquad\text{ with $0 \leq j_k\leq l$ and $\sum_{k=1}^N (j_k+1)=s$,}
$$
where the covariant indexes are all completely contracted with
$g^{ij}$.
\end{dfnz}

\begin{rem}
See the paper~\cite[Sect.~2]{mant5} for more details on these
polynomials and the geometric interpretation of the rescaling
  order.
Notice that, differently from~\cite{mant5}, here we need to consider in 
$\pol^s(\nabla^l\BBB)$ monomials of two types, because of the
codimension higher than one.
\end{rem}

We advise the reader that the polynomials $\pol^s$ and  
$\qol^s$ may vary from line to line, 
and similarly the constants (usually indicated by $C$).

\section{Evolving Plane Curves}\label{seccurves}

Let $\gamma\in C^\infty(S^1;\R^2)$ be a regular immersed closed 
curve in the plane $\R^2$. 
Let $\tau={\gamma_x}/{\vert \gamma_x\vert}=\gamma_s$ and
$\nu={\mathrm R}\tau$ be respectively the tangent and the normal to the
curve $\gamma$, where ${\mathrm R}$ is the counterclockwise
rotation of $\pi/2$ in the plane, and $\gamma_x = \partial_x \gamma$.\\
We recall that $\partial_s ={\partial_x}/{\vert \gamma_x\vert}$ and
\begin{equation}\label{freser}
\partial_s \tau = \kappa \nu, \qquad \partial_s 
\nu = - \kappa \tau
\end{equation}
where $\k$ is the curvature of $\gamma$.
In the sequel we let $\L=\L(\gamma)=\int_{\gamma}1\,ds$ 
be the length of the curve $\gamma$. 

Let us consider the functional
\begin{equation*}
{\mathcal G}^\eps(\gamma) = \int_{\gamma} \left(1 + \eps 
\kappa^2\right)\, ds\,, 
\end{equation*}
which is obtained from ${\mathcal G}^\eps_3$ (with $n=m=1$)
by replacing $3\eps$ with $\eps$. 
Set 
\begin{equation*}
\Eul = -\kappa +2 \eps \partial_s^2 \kappa + \eps \kappa^3\,. 
\end{equation*}
Then the gradient flow by ${\mathcal G}^\eps$ is given by 
a smooth map $\gamma: S^1 \times [0,+\infty)\to \R^2$ which is an 
immersion for any $t \in [0,+\infty)$, equals a given immersion $\gamma_0$
at time $t=0$, and satisfies
\begin{equation}\label{eqevoleps}
\partial_t \gamma =
-\Eul ~\nu\,,
\end{equation}
where $\partial_t = \frac{\partial}{\partial t}$. For notational 
simplicity, we omit the dependence of $\gamma$ on $\eps$.
\begin{lemma}\label{lemzero}
We have
\begin{equation*}
\partial_s \partial_t \gamma = 
-(\partial_s\Eul) \nu + \kappa \Eul \tau\,,
\end{equation*}
in particular 
\begin{equation}\label{inpa}
\langle \partial_s \partial_t \gamma, 
\tau\rangle =\kappa \Eul\,.
\end{equation}
\end{lemma}
\begin{proof}
It follows from equations~\eqref{freser} and the evolution
equation~\eqref{eqevoleps}.
\end{proof}

\begin{lemma}\label{lemsLkappa}
Let $\gamma$ be a smooth closed curve, then
\begin{equation}\label{dislkappa}
\frac{1}{\L} \leq \frac{1}{4\pi^2} 
\int_\gamma \kappa^2\,ds\,.
\end{equation}
\end{lemma}
\begin{proof}
By Borsuk and Schwartz--H\"older inequalities we have 
$$
2\pi \leq \int_\gamma \vert \kappa\vert ~ds \leq \left( \int_\gamma 
\kappa^2\,ds\right)^{1/2} \L^{1/2}\,.
$$
\end{proof}

\begin{lemma}\label{lemscambio}
The following commutation rule holds:
\begin{equation}\label{eqscambio}
\partial_t \partial_s = \partial_s \partial_t 
-\kappa \Eul\partial_s\,.
\end{equation}
\end{lemma}
\begin{proof}
Observing that 
$\frac{\partial_t \partial_x}{\vert \gamma_x\vert}=
\frac{\partial_x}{\vert \gamma_x\vert} \partial_t = \partial_s \partial_t$, 
we have
\begin{align*}
\partial_t \partial_s 
=&\,\partial_t \left( \frac{\partial_x}{\vert \gamma_x\vert}\right)
= \frac{\partial_t \partial_x}{\vert \gamma_x\vert} 
-\frac{\langle\gamma_x, \partial_t \gamma_x\rangle \partial_x}{\vert 
\gamma_x\vert^3}
=\frac{\partial_x}{\vert \gamma_x\vert} \partial_t 
-\left\langle\frac{\gamma_x}{\vert \gamma_x\vert},
\frac{\partial_t \gamma_x}{\vert \gamma_x\vert}\right\rangle
\frac{\partial_x}{\vert \gamma_x\vert}\\
= &\,\partial_s \partial_t - 
\langle \tau, \partial_s \partial_t \gamma\rangle \partial_s\,.
\end{align*}
Then the commutation rule~\eqref{eqscambio} follows from equation~\eqref{inpa}. 
\end{proof}

\begin{lemma}\label{lemett}
We have 
\begin{equation}\label{dertkappa}
\partial_t \kappa = -\partial^2_s \Eul 
-\kappa^2 \Eul = \partial^2_s \k +\k^3-2\eps \partial_s^4 \k 
-6\eps\k (\partial_s \k)^2-5\eps\k^2 \partial^2_s \k-\eps\k^5\,.
\end{equation}
\end{lemma}
\begin{proof}
We have
\begin{equation*}
\partial_t \kappa = \partial_t \langle \partial_s \tau, \nu \rangle
= \langle \partial_t \partial_s \tau, \nu \rangle\,.
\end{equation*}
Therefore, using formula~\eqref{eqscambio} we have
\begin{equation*}
\partial_t \kappa 
= \langle \partial_s  \partial_t \partial_s \gamma, \nu \rangle 
-\kappa \Eul  \langle 
\partial_s \tau, \nu \rangle
= \langle \partial^2_s \partial_t \gamma, \nu\rangle -
\langle \partial_s [\kappa \Eul \partial_s \gamma],\nu \rangle - 
\kappa^2 \Eul\,.
\end{equation*}
Using the evolution law~\eqref{eqevoleps} we get
$$
\langle \partial^2_s \partial_t \gamma, \nu\rangle  = 
- \langle \partial^2_s (\Eul 
\nu), \nu\rangle = - \partial^2_s \Eul 
+\Eul \langle \partial_s (\kappa \tau), \nu\rangle = 
-\partial^2_s \Eul + \kappa^2 \Eul\,.
$$
In addition,
$$
\langle \partial_s [\kappa \Eul \partial_s \gamma],\nu \rangle =
\kappa \Eul \langle \partial_s \tau, \nu \rangle =
\kappa^2 \Eul\,. 
$$
Hence
$\partial_t \kappa = - 
\partial^2_s \Eul - \kappa^2\Eul$ and the last equality
in~\eqref{dertkappa} follows by expanding $\Eul$.
\end{proof}

\begin{rem}\label{primocontrolloapost}
For $\eps=0$, formula~\eqref{dertkappa} gives
the well known evolution equation 
$\kappa_t = \partial_s^2 \kappa + \kappa^3$, valid for motion by 
curvature,
see~\cite[Lemma 3.1.6]{gaha1}.
\end{rem}

We recall now the following interpolation inequalities for closed
curves, see~\cite[pag.~93]{aubin0}.
\begin{prop}\label{propinter}
Let $\gamma$ be a regular closed curve in
  $\R^2$ with finite length $\L$. Let
$u$ be a smooth function 
defined on
  $\gamma$, $m\geq1$ and $p\in[2,+\infty]$. 
If $n\in\{0,\dots, m-1\}$ we have the estimates
\begin{equation}\label{int1}
    {\Vert\partial_s^n u\Vert}_{L^p}
    \leq C_{n,m,p}
      {\Vert\partial_s^m  u\Vert}_{L^2}^{\sigma}
      {\Vert u\Vert}_{L^2}^{1-\sigma}+
      \frac{B_{n,m,p}}{\L^{m\sigma}}{\Vert u\Vert}_{L^2}\,,
  \end{equation}
where
\begin{equation*}
\sigma = \frac{n+1/2-1/p}{m}\in[0,1)
\end{equation*}
and the constants $C_{n,m,p}$ and $B_{n,m,p}$ are independent of $\gamma$.
\end{prop}
Clearly inequalities~\eqref{int1} hold with uniform constants if
applied to a family of curves having lengths uniformly 
bounded below by some positive value.

\begin{rem}
In the special case $p=+\infty$,
we have $\sigma=\frac{n+1/2}{m}$, and
\begin{equation*}
    {\Vert\partial_s^n u\Vert}_{L^\infty}
    \leq C_{n,m}
      {\Vert\partial_s^m  u\Vert}_{L^2}^{\sigma}
      {\Vert u\Vert}_{L^2}^{1-\sigma}+
      \frac{B_{n,m}}{{\mathrm L}^{m\sigma}}{\Vert
        u\Vert}_{L^2}\,.
\end{equation*}
\end{rem}

\begin{rem}\label{remarema}
In the particular case $n=0$, $m=2$, $p=6$ we get 
$\sigma=1/6$ and 
\begin{equation*}
    {\Vert u\Vert}_{L^6}
    \leq C 
      {\Vert \partial^2_s u\Vert}_{L^2}^{\frac{1}{6}}
      {\Vert u\Vert}_{L^2}^{\frac{5}{6}}+
      \frac{C}{{\mathrm L}^{\frac{1}{3}}}{\Vert u\Vert}_{L^2}\,,
\end{equation*}
for some $C>0$, hence, by means of Young inequality 
$\vert xy\vert  \leq \frac{1}{\alpha}\vert x\vert^a + \frac{1}{b} 
\vert y\vert^b$, $1/a+1/b=1$, choosing $a=b=2$, $x=\sqrt{2}\Vert 
\partial_s^2 u\Vert_{L^2}^{1/2}$ and 
$y= \frac{1}{\sqrt{2}}\Vert u\Vert_{L^2}^{5/2}$, we obtain
\begin{equation}\label{int4}
    \int_\gamma u^6\, ds \le \int_\gamma (\partial_s^2 u)^2\, ds + 
    C \left(\int_\gamma u^2\, ds\right)^5
+     \frac{C}{{\mathrm L}^2}\left(\int_\gamma u^2\, ds\right)^3\,.
\end{equation}
In the particular case $n=0$, $m=1$, $p=4$ we get 
$\sigma=1/4$ and 
\begin{equation*}
    {\Vert u\Vert}_{L^4}
    \leq C 
      {\Vert \partial_s u\Vert}_{L^2}^\frac{1}{4}
      {\Vert u\Vert}_{L^2}^\frac{3}{4}+
      \frac{C}{{\mathrm L}^{\frac{1}{4}}}{\Vert u\Vert}_{L^2}\,,
\end{equation*}
hence, reasoning as before,
\begin{equation}\label{int6}
    \int_\gamma u^4\, ds \le \int_\gamma (\partial_s u)^2\, ds + 
    C \left(\int_\gamma u^2\, ds\right)^{3}
+     \frac{C}{{\mathrm L}}\left(\int_\gamma u^2\, ds\right)^2\,.
\end{equation}
\end{rem}

We are now ready for the estimates. We recall that 
\begin{equation}\label{eqevollung}
\partial_t\,ds =
\kappa \Eul  \,ds= 
(-\kappa^2+2\eps\kappa \partial_s^2 \kappa+\eps\kappa^4)\,ds.
\end{equation}
\begin{lemma}\label{lemintkapqua}
We have
\begin{equation}\label{evollduenormkappa}
\partial_t \int_{\gamma} \kappa^2 \,ds = 
\int_{\gamma} \big(- 2(\partial_s \k)^2 
+\k^4-4\eps (\partial^2_s \k)
-\eps\k^6 
-4\eps\k^3\partial_s^2 \k\big) \,ds\,.
\end{equation}
\end{lemma}
\begin{proof}
From~\eqref{eqevollung} and Lemma~\ref{lemett} we get
\begin{align*}
\partial_t \int_{\gamma} \kappa^2 \,ds 
=&\, 2\int_{\gamma} \kappa \partial_t \kappa \,ds + \int_{\gamma}
(-\kappa^4+2\eps\kappa^3 \partial_s^2 \kappa+\eps\kappa^6)\,ds\\
=&\, 2\int_{\gamma} \kappa\big(\partial_s^2 \k+\k^3-2\eps\partial_s^4 \k 
-6\eps\k (\partial_s \k)^2-5\eps\k^2 \partial^2_s \k-\eps\k^5\big)\,ds\\
&\,+ \int_{\gamma}(-\kappa^4+2\eps\kappa^3\partial^2_s \kappa
+\eps\kappa^6)\,ds\\
=&\, \int_{\gamma} \big(2\kappa \partial_s^2 \k+\k^4-4\eps\k 
\partial_s^4 \k 
-12\eps\k^2 (\partial_s \k)^2-8\eps\k^3 \partial^2_s \k-\eps\k^6\big)\,ds.
\end{align*}
Therefore, integrating by parts, we obtain
\begin{align*}
\partial_t \int_{\gamma} \kappa^2 \,ds
=&\, \int_{\gamma} \big(-2 (\partial_s \kappa)^2 +\k^4
-4\eps (\partial^2_s 
\kappa)^2
-\eps\kappa^6-12\eps\k^2 (\partial_s \k)^2-8\eps\k^3 \partial_s^2 
\k\big)\,ds\\
=&\, \int_{\gamma} \big(-2 (\partial_s \kappa)^2 
+\k^4
-4\eps(\partial_s^2\kappa)^2
-\eps\kappa^6-4\eps\k^3 \partial_s^2\k\big)\,ds,
\end{align*}
where in the last equality we used the fact that 
$-3\int_\gamma\k^2 
(\partial_s \k)^2\,ds=\int_\gamma\k^3\partial_s^2\k\,ds$.
\end{proof}

\begin{prop}\label{ltwoboundkappa}
The following estimate holds
\begin{equation}
\label{ineqltwonorm}
\partial_t \int_{\gamma}
\kappa^2 
\,ds
\le 
C \left( \int_{\gamma}  \kappa^2\,ds
\right)^{3}
+ C \left( \int_{\gamma}  \kappa^2\,ds
\right)^5\,,
\end{equation}
where $C$ is a constant independent of $\varepsilon$.
\end{prop}
\begin{proof}
Adding to the right hand side of equation~\eqref{evollduenormkappa} the
positive quantity $2\eps(\partial_s^2\kappa + \kappa^3)^2$ 
we get
$$
\partial_t \int_{\gamma} \kappa^2 \,ds \leq 
\int_{\gamma} \big(- 2 (\partial_s \k)^2 +\k^4 
-2\eps (\partial_s^2\k)^2 
+ 
\eps\k^6\big) 
\,ds\,.
$$
Using now inequalities~\eqref{int4} and~\eqref{int6} we obtain
\begin{align*}
\partial_t \int_{\gamma} \kappa^2 
\leq &\,\int_{\gamma} \big(- (\partial_s 
\k)^2
-\eps (\partial_s^2 \k)^2 
\big) \,ds
+C \eps \left(\int_{\gamma} \kappa^2\,ds
\right)^5
+ \frac{C \eps}{\L^2} \left( \int_{\gamma} \kappa^2\,ds
\right)^3
\\
& + C \left( \int_{\gamma} \kappa^2\,ds
\right)^{3}
+ \frac{C}{\L} \left( \int_{\gamma} \kappa^2\,ds
\right)^2
\\
\leq &\,C \left(\int_{\gamma} \kappa^2\,ds
\right)^5
+C\left( \int_{\gamma} \kappa^2\,ds
\right)^3
+ \frac{C}{\L} \left( \int_{\gamma} \kappa^2\,ds
\right)^2
+ \frac{C}{\L^2}\left( \int_{\gamma} \kappa^2\,ds
\right)^3\\
\leq &\,C \left(\int_{\gamma} \kappa^2\,ds
\right)^5
+C\left( \int_{\gamma} \kappa^2\,ds
\right)^3\,,
\end{align*}
where we supposed $\eps<1$ and in the last inequality we used 
the geometric 
estimate~\eqref{dislkappa}.
\end{proof}

We deal now with the higher derivatives of the curvature.\\
Since here we are working in dimension and codimension one, for 
the rest of
this section all polynomials in the curvature $\k$ and its derivatives 
are completely contracted, that is they belong to 
the ``family'' $\qol^r(\partial^l_s\k)$ (see Definition~\ref{polqol});
moreover, every of their monomials is of the form
$$
\prod_{i=1}^N\partial_s^{j_i}\k
\qquad\text{ with $0 \leq j_i\leq l$ and $N\geq1$}
$$
with
$$
r=\sum_{i=1}^N (j_i+1)\,,
$$
as the $*$ product in this case is simply the usual product.

\begin{lemma}\label{kexpr}
For any $j\in\NN$ the following formula holds: 
\begin{equation}\label{kappas}
\dert\ders^j\k=\ders^{j+2}\k + \qol^{j+3}(\ders^{j}\k)
-2\eps \ders^{j+4}\k -5\eps\k^2\ders^{j+2}\k
+ \eps \qol^{j+5}(\ders^{j+1}\k)\,.
\end{equation}
\end{lemma}
\begin{proof} 
We argue by induction on $j$.\\ 
The case $j=0$ in~\eqref{kappas} is
equation~\eqref{dertkappa}, where 
$\qol^5(\partial_s \kappa) = - 6 \kappa (\partial_s \kappa)^2 - 
\kappa^5$.\\
Suppose that~\eqref{kappas} holds for
$(j-1)$; using the commutation rule~\eqref{eqscambio} we get 
\begin{align*}
\dert\ders^j \k=&\, \ders\dert\ders^{j-1}\k 
+ \kappa (\kappa - 2 \eps \partial_s^2\kappa - \eps \kappa^3) \ders^j 
\k\\
=&\, \ders \left[\ders^{j+1}\k + \qol^{j+2}(\ders^{j-1}\k)
-2\eps\ders^{j+3}\k -5\eps\k^2\ders^{j+1}\k 
+ \eps \qol^{j+4}(\ders^{j}\k)\right]\\
& + \qol^{j+3}(\ders^j\k)+\eps\qol^{j+5}(\ders^{j}\k)\,,
\end{align*}
where we expressed $\qol^{j+3}(\partial_s^j \kappa) = \kappa^2 
\partial_s^j \kappa$ and $\qol^{j+5}(\partial_s^j \kappa) = -
(2\kappa \partial_s^2 \kappa + \kappa^4)\partial_s^j \kappa$. Hence, 
we deduce
\begin{equation*}
\dert\ders^j \k
= \ders^{j+2}\k + \qol^{j+3}(\ders^{j}\k)
-2\eps \ders^{j+4}\k 
-5\eps\k^2\ders^{j+2}\k
+ \eps \qol^{j+5}(\ders^{j+1}\k)\,, 
\end{equation*}
which gives the inductive step.
\end{proof}

\begin{lemma}\label{lemevolnormgene}
For any $j\in\NN$ we have
\begin{align*}
\dert \int_{\gamma} |\ders^j \k|^2 \,ds 
= & -2\int_{\gamma} \vert\ders^{j+1} \k\vert^2\,ds
-4\eps \int_{\gamma} \vert\ders^{j+2} \k\vert^2 \,ds
\\
& +
 \int_{\gamma} \qol^{2j+4}(\ders^j\k)\,ds 
+ \eps \int_{\gamma} \qol^{2j+6}(\ders^{j+1}\k)\,ds\,.
\end{align*}
\end{lemma}
\begin{proof}
Using~\eqref{eqevollung},~\eqref{kappas} and integrating 
by parts we deduce
\begin{align}\label{evolint000}
\dert \int_{\gamma} |\ders^j \k|^2 \,ds 
=&\,2\int_{\gamma} \ders^j \k\, \dert \ders^j \k\,ds 
+ \int_{\gamma} |\ders^j \k|^2 \kappa \Eul \,ds\\
\nonumber 
=&\, 2\int_{\gamma} \ders^j \k\, \big(\ders^{j+2}\k 
+\qol^{j+3}(\ders^j\k)\bigr)\,ds\\
&\,+\varepsilon\int_{\gamma}2\ders^j \k\, 
\big(-2\ders^{j+4}\k -5\k^2\ders^{j+2}\k
+ \qol^{j+5}(\ders^{j+1}\k)\big)\,ds\nonumber\\
\nonumber 
&\,- \int_{\gamma} |\ders^j \k|^2\k(\k-2\eps 
\partial_s^2 \k-\eps\k^3)\,ds\\
\nonumber 
=&\,-2\int_{\gamma} \left(\vert\ders^{j+1} \k\vert^2 
+ \qol^{2j+4}(\ders^j\k)\right)\,ds\\
&\,-4\eps \int_{\gamma} \left(\vert\ders^{j+2} \k\vert^2 
+ \qol^{2j+6}(\ders^{j+1}\k)\right)\,ds\,.\nonumber
\end{align}
\end{proof}

\begin{prop}\label{teobello} For any $j\in\NN$ we have the
  $\eps$--independent estimate, for $\eps<1$,
\begin{equation}\label{eps1est}
\dert \int_{\gamma} |\ders^j \k|^2 \,ds
\le C \left(\int_{\gamma}\k^2\,ds\right)^{2j+3} 
+ C \left(\int_{\gamma}\k^2\,ds\right)^{2j+5} + C
\end{equation}
where the constant $C$ depends only on $1/\L$.
\end{prop}
\begin{proof}
We estimate the term $\int_\gamma\qol^{2j+4}(\ders^j\k)\,ds$ as
in~\cite[Sect.~7]{mant5}.
By definition, we have
$$
\qol^{2j+4}(\ders^j\k)=\sum_m \prod_{l=1}^{N_m}\ders^{c_{ml}}\k
$$
with all the $c_{ml}$ less than or equal to $j$ and 
$$
\sum_{l=1}^{N_m} (c_{ml}+1)=2j+4
$$
for every $m$. Hence,
$$
\vert\qol^{2j+4}(\ders^j\k)\vert\leq\sum_m \prod_{l=1}^{N_m}\vert\ders^{c_{ml}}\k\vert
$$
and setting
$$
Q_m =\prod_{l=1}^{N_m}\vert\ders^{c_{ml}}\k\vert\,,
$$
we clearly obtain
$$
\int_\gamma \vert\qol^{2j+4}(\ders^j\k)\vert\,ds\leq\sum_m\int_\gamma Q_m\,ds\,.
$$
We now estimate any term $Q_m$ via interpolation inequalities.
After collecting derivatives of the same order in $Q_m$ we can write
\begin{equation}\label{condos}
Q_m = \prod_{i=0}^{j}\vert\ders^i \k\vert^{\alpha_{mi}} \qquad
\text{ with } \quad \sum_{i=0}^{j}\alpha_{mi}(i+1) = 2j+4\,.
\end{equation}
Then
\begin{equation*}
\int_{\gamma} Q_m\,ds=
\int_{\gamma}\prod_{i=0}^{j}\vert\ders^i
\k\vert^{\alpha_{mi}}\,ds
\leq
\prod_{i=0}^{j}\left(\int_{\gamma}\vert\ders^i
    \k\vert^{\alpha_{mi}\lambda_i}\,ds\right)^{\frac{1}{\lambda_i}}
=
\prod_{i=0}^{j}\Vert\ders^i
\k\Vert_{L^{\alpha_{mi}\lambda_i}}^{\alpha_{mi}}
\end{equation*}
where the values $\lambda_i$ are chosen as follows: $\lambda_i=0$ if 
$\alpha_{ji}=0$ (in this case the corresponding term is not present in 
the product) and
$\lambda_i=\frac{2j+4}{\alpha_{mi}(i+1)}$ if
$\alpha_{mi}\not=0$. Clearly,
$\alpha_{mi}\lambda_i=\frac{2j+4}{i+1}\geq\frac{2j+4}{j+1}>2$ and by
the condition in~\eqref{condos},
$\sum_{\genfrac{}{}{0pt}{}{i=0}{\lambda_i\not=0}}^j
\frac{1}{\lambda_i}=\sum_{\genfrac{}{}{0pt}{}{i=0}{\lambda_i\not=0}}^j
\frac{\alpha_{mi}(i+1)}{2j+4}=1$.\\
As $\alpha_{mi}\lambda_i>2$ these values are allowed as exponents
$p$ in inequality~\eqref{int1} and taking $m = j+1$, 
$n=i$, $u=\k$, we get
$$
    {\Vert\ders^{i} \k\Vert}_{L^{\alpha_{mi}\lambda_i}}
\leq C{\Vert\ders^{j+1}\k\Vert}_{L^2}^{\sigma_{mi}}
{\Vert\k\Vert}_{L^2}^{1-\sigma_{mi}}+\frac{C}{L^{(j+1)\sigma_{mi}}}\Vert\k\Vert_{L^2}
\leq C\left(\Vert\ders^{j+1}\k\Vert_{L^2}+\Vert\k\Vert_{L^2}\right)^{\sigma_{mi}}
{\Vert\k\Vert}_{L^2}^{1-\sigma_{mi}}
$$
with
\begin{equation*}
\sigma_{mi}=\frac{i+1/2-1/(\alpha_{mi}\lambda_i)}{j+1}
\end{equation*}
and the constant $C$ depends only on $1/\L$.

Multiplying together all the estimates,
\begin{align}
\int_{\gamma} Q_m\,ds
\leq
\,&\,C\prod_{i=0}^{j}
\left(\Vert\ders^{j+1}\k\Vert_{L^2}+\Vert\k\Vert_{L^2}\right)^{\alpha_{mi}
\sigma_{mi}}
{\Vert\k\Vert}_{L^2}^{\alpha_{mi}(1-\sigma_{mi})}\label{equ1000}\\
\leq\,&\,C\left(\Vert\ders^{j+1}\k\Vert_{L^2}
+\Vert\k\Vert_{L^2}\right)^{\sum_{i=0}^j\alpha_{mi}\sigma_{mi}}
{\Vert\k\Vert}_{L^2}^{\sum_{i=0}^j\alpha_{mi}(1-\sigma_{mi})}\,.\nonumber
\end{align}
Then we compute
$$
\sum_{i=0}^j\alpha_{mi}\sigma_{mi}=\sum_{i=0}^j
\frac{\alpha_{mi}(i+1/2)-1/\lambda_i}{j+1}=
\frac{\sum_{i=0}^j\alpha_{mi}(i+1/2)-1}{j+1}
$$
and using again the rescaling condition in~\eqref{condos},
$$
\sum_{i=0}^j\alpha_{mi}\sigma_{mi} 
=\frac{4j+6-\sum_{i=0}^j\alpha_{mi}}{2(j+1)}\,.
$$
Since
$$
\sum_{i=0}^j\alpha_{mi}\geq\sum_{i=0}^j\alpha_{mi}\frac{i+1}{j+1}=\frac{2j+4}{j+1}
$$
we get
$$
\sum_{i=0}^j\alpha_{mi}\sigma_{mi}\leq\frac{2j^2+4j+1}{(j+1)^2}=2-\frac{1}{(j+1)^2}<2\,.
$$
Hence, we can apply Young inequality to the product in the last term
of inequality~\eqref{equ1000}, in order to get the exponent 2 on the
first quantity, that is,
$$
\int_{\gamma} Q_m\,ds
\leq \frac{\delta_m}{2}\left(\Vert\ders^{j+1}\k\Vert_{L^2}
+\Vert\k\Vert_{L^2}\right)^{2}+
{\Vert\k\Vert}_{L^2}^\beta
\leq \delta_m\int_{\gamma} \vert\ders^{j+1}\k\vert^2\,ds+
\delta_m\int_{\gamma} \vert\k\vert^2\,ds+
{\Vert\k\Vert}_{L^2}^\beta\,,
$$
for arbitrarily small $\delta_m>0$ and where $\beta$ is given by
\begin{align*}
\beta=
&\,\sum_{i=0}^j\alpha_{mi}(1 
-\sigma_{mi})\frac{1}{1-\frac{\sum_{i=0}^j\alpha_{mi}\sigma_{mi}}{2}}\\
=&\,\frac{2\sum_{i=0}^j\alpha_{mi}(1-\sigma_{mi})}{2-\sum_{i=0}^j\alpha_{mi}\sigma_{mi}}\\
=&\,\frac{2\sum_{i=0}^j\alpha_{mi}-\frac{4j+6-\sum_{i=0}^j\alpha_{mi}}{j+1}}
{2-\frac{4j+6-\sum_{i=0}^j\alpha_{mi}}{2(j+1)}}\\
=&\,2\frac{2(j+1)\sum_{i=0}^j\alpha_{mi}-4j-6+\sum_{i=0}^j\alpha_{mi}}
{4j+4-4j-6+\sum_{i=0}^j\alpha_{mi}}\\
=&\,2\frac{(2j+3)\sum_{i=0}^j\alpha_{mi}-2(2j+3)}
{\sum_{i=0}^j\alpha_{mi}-2}\\
=&\,2(2j+3)\,.
\end{align*}
Therefore we conclude
$$
\int_{\gamma} Q_m\,ds
\leq \delta_m\int_{\gamma} \vert\ders^{j+1}\k\vert^2\,ds+
\delta_m\int_{\gamma} \k^2\,ds+ C\left(\int_{\gamma}
  \k^2\,ds\right)^{2j+3}\,.
$$
Repeating this argument for all
the $Q_m$ and choosing suitable $\delta_m$ whose sum over $m$ is less than
one, we conclude that there exists a constant $C$ depending only on
$1/\L$ and $j\in\NN$ such that
\begin{equation*}
\int_{\gamma}\qol^{2j+4}(\ders^j\k)\,ds\,\leq
\,\int_\gamma\vert\ders^{j+1}\k\vert^2\,ds+
C\left(\int_\gamma\k^2\right)^{2j+3}+C\,.
\end{equation*}

Reasoning similarly for the term $\qol^{2j+6}(\ders^{j+1}\k)$, we obtain
\begin{equation*}
\int_{\gamma} \qol^{2j+6}(\ders^j\k)\,ds \le 
\int_{\gamma} \vert\ders^{j+2} \k\vert^2 \,ds
+ C\left(\int_{\gamma}\k^2\,ds\right)^{2j+5} + C\,.
\end{equation*}
Hence, from~\eqref{evolint000} we get
\begin{align*}
\dert \int_{\gamma} |\ders^j \k|^2 \,ds 
\le &\,-\int_{\gamma} \vert\ders^{j+1} \k\vert^2\,ds
-\eps \int_{\gamma} \vert\ders^{j+2} \k\vert^2\,ds\\
& +C\left(\int_{\gamma}\k^2\,ds\right)^{2j+3} 
+ C\eps \left(\int_{\gamma}\k^2\,ds\right)^{2j+5} + C\\
\le &\,C\left(\int_{\gamma}\k^2\,ds\right)^{2j+3} 
+ C\left(\int_{\gamma}\k^2\,ds\right)^{2j+5} + C
\end{align*}
when $\eps<1$ and the constant $C$ depends only on $1/\L$.
\end{proof}

By means of Propositions~\ref{ltwoboundkappa} and~\ref{teobello} we
have then the following result.

\begin{teo}\label{teoteo1} For any $j\in\NN$ there exists a smooth
  function $\ZZZ^j:\R\to(0,+\infty)$ such that
\begin{equation*}
\dert \int_{\gamma} |\ders^j \k|^2 \,ds 
\leq \ZZZ^j\left(\int_{\gamma}\k^2 \,ds\right)
\end{equation*}
for every $\eps<1$ and curve $\gamma$ evolving by the gradient of the
functional ${\mathcal G}^\eps$.
\end{teo}
\begin{proof}
The statement clearly follows by Propositions~\ref{ltwoboundkappa}
and~\ref{teobello}, since by Lemma~\ref{lemsLkappa} the quantity
$1/\L$ is controlled by $\int_\gamma\k^2\,ds$.\\
The smoothness of the functions $\ZZZ^j$ is obtained choosing possibly
slightly larger constants in inequalities~\eqref{eps1est}
and~\eqref{ineqltwonorm}.
\end{proof}

This proposition, like its analogue for the general 
case (Theorem~\ref{teoteo3}), is the key tool in order to get
$\eps$--independent compactness estimates. Indeed, for example, 
one can see that, by an ODE's argument, since all the flows (letting
$0<\eps<1$ vary) start from a common initial smooth curve, fixing any
$j\in\NN$, there exists a common positive interval of time such that 
all the quantities $\Vert\ders^i\k\Vert_{L^2}$, for
$i\in\{0,\dots,j\}$ are equibounded. This will allow us to get
compactness and $C^\infty$ convergence to the mean curvature flow as
$\eps\to0$.

\section{The General Case}\label{secgen}

If $k>[n/2]+2$ it is shown in~\cite{mantemin1} 
that for every $\varepsilon>0$ all the flows $\varphi^\varepsilon$,
associated with the functionals
$$
\DG^\varepsilon(\varphi) = \int_M
\left(1+\varepsilon\vert A^k\vert^2\right)\,d\mu\,,
$$
and starting from a common initial $n$--dimensional smooth 
compact immersed submanifold, 
are smooth for every positive time.

By means of Theorem~4.5 and Theorem~5.9 in~\cite{ambman1} and the
results of~\cite{mantemin1}, the {\em first variation} of the functional 
$\DG^\varepsilon$ is given by
\begin{equation*}
\EEE^\varepsilon=-\HHH+ 2\varepsilon
k(-1)^{k-1}\Bigl(\overset{\text{$( k -2)$--times}}{\overbrace{\Delta^{M}\comp
\Delta^{M}\comp\,\dots\,\comp\Delta^{M}}}\,\,\HHH\Bigr)^\perp
+\varepsilon\qol^{2k-3}(\nabla^{2k-5}\BBB)^\perp
\end{equation*}
where $\qol^{2k-3}(\BBB)$ takes values in $\R^{n+m}.$\\
Here we denote with $\Delta^{M}$ the Laplacian of the smooth 
compact Riemannian
manifold without boundary $M=(M,g)$, where $g$ is the metric induced on
$M$ by the immersion.

Then we have a  solution of the geometric evolution problem for any initial
smooth immersion $\varphi_0:M\to\R^{n+m}$, that is, a 
smooth
function $\varphi^\varepsilon:M\times[0,+\infty)\to\R^{n+m}$ such that
\begin{enumerate}
\item the map $\varphi^\varepsilon(\cdot,t):M\to\R^{n+m}$ is an
  immersion for every $t\in[0,+\infty)$;
\item $\varphi^\varepsilon(p,0)=\varphi_0(p)$ for every $p\in M$;
\item the following parabolic system is satisfied
\begin{equation*}
\frac{\partial\varphi^\varepsilon}{\partial t}=-\EEE^\varepsilon=
\HHH+ 2\varepsilon k(-1)^{k}
\Bigl(\overset{\text{$( k -2)$--times}}{\overbrace{\Delta^{M}\comp
\Delta^{M}\comp\,\dots\,\comp\Delta^{M}}}\,\,\HHH\Bigr)^\perp
+\varepsilon\qol^{2k-3}(\nabla^{2k-5}\BBB)^\perp\,.
\end{equation*}
\end{enumerate}

\section{Evolution of Geometric Quantities}\label{evolgeosec}

We work out some evolution equations for the
geometric quantities under the flow by the gradient of
$\DG^\varepsilon$.

In general, if a family of immersed manifolds
$\varphi(\cdot, t):M\to\R^{n+m}$ moves by $\partial_t\varphi=-\EEE$, 
with
the field $\EEE$ normal, we have
\begin{equation*}
 \derpar t {g_{ij} } = 2  \left\langle \BBB_{ij} ,
   \EEE\right\rangle\,,\qquad
 \derpar t {g^{ij}} = 
-2 g^{is} \left\langle \BBB_{sl} , \EEE \right\rangle g^{lj}\,.
\end{equation*}
Now for the Christoffel symbols $\Gamma_{ij}^s = 
\left\langle \frac{\partial^2 {\varphi}}{\partial x_i \partial x_j},\derpar
  x_l {\varphi} \right\rangle g^{ls} $ we have
\begin{equation*}
 \derpar t {\Gamma_{ij}^s} = -
 \left\langle {\frac{\partial^2 \EEE}{\partial x_i \partial x_j}},
  \derpar {x_l} {\varphi} \right\rangle g^{ls} -
 \left\langle \frac{\partial^2 {\varphi}}{\partial x_i \partial x_j},\derpar x_l {\EEE} \right\rangle g^{ls} +
 \left\langle \frac{\partial^2 {\varphi}}{\partial x_i  \partial 
x_j},\derpar x_l {{\varphi}} \right\rangle \derpar t {g^{ls}}\,.
 \end{equation*}
Then, supposing to work in normal coordinates, $\frac{\partial^2
  {\varphi}}{\partial x_i \partial x_j}=\BBB_{ij}$ is a normal vector, hence
\begin{align*}
 \derpar t {\Gamma_{ij}^s} = 
&\, -\left\langle {\frac{\partial^2 \EEE}{\partial x_i \partial x_j}},
  \derpar {x_l} {\varphi} \right\rangle g^{ls} -
 \left\langle \BBB_{ij},\derpar {x_l} {\EEE}
\right\rangle g^{ls}\\
= &\,- {\frac{\partial^2 \,\,}{\partial x_i \partial x_j}}\left\langle
      \EEE, \frac{\partial\varphi}{\partial x_l} \right\rangle g^{ls}
    +\left\langle\frac{\partial\EEE}{\partial x_i},
      \frac{\partial^2\varphi}{\partial x_j\partial x_l} \right\rangle
    g^{ls}
+\left\langle\frac{\partial\EEE}{\partial x_j},
      \frac{\partial^2\varphi}{\partial x_i\partial x_l} \right\rangle
    g^{ls}\\
&\,+\left\langle
      \EEE, \frac{\partial^3\varphi}{\partial x_i\partial x_j\partial
        x_l} \right\rangle g^{ls}
- \left\langle \BBB_{ij},\frac{\partial\EEE}{\partial x_l}
\right\rangle g^{ls}\\
= & \, 
\left\langle \BBB_{jl},\nabla_i{\EEE}\right\rangle g^{ls}
+
\left\langle \BBB_{il},\nabla_j {\EEE}\right\rangle g^{ls}
+
\left\langle      \EEE, \nabla_i\BBB_{jl} \right\rangle g^{ls}
- 
\left\langle \BBB_{ij},\nabla_l{\EEE}\right\rangle g^{ls}\,.
\end{align*}

\begin{rem}\label{rem3bbb}
By this computation, since the Christoffel
symbols are symmetric in the ${ij}$--indices, the covariant 3--tensor
$(\nabla\BBB)^\perp$ is symmetric (as in the codimension one case).
\end{rem}

Then, we compute the evolution of $\BBB$,
\begin{align}\label{enzo}          
\derpar t {\BBB_{ij}} = 
&\, \frac{\partial\,}{\partial t}\left(\frac{\partial^2 
{\varphi}}{\partial x_i \partial x_j}
-\Gamma_{ij}^s \derpar {x_s} {\varphi} \right) 
=  
-\frac{\partial^2 \EEE}{\partial x_i \partial x_j} 
 - \derpar t {\Gamma_{ij}^s} \derpar {x_s} {\varphi}
+ \Gamma_{ij}^s \derpar {x_s} {\EEE} 
\\
 = &\, -\nabla^2_{ij} \EEE - \derpar t {\Gamma_{ij}^s} 
\derpar {x_s} {\varphi} 
= - \nabla^2_{ij} \EEE
+\fol_s(\nabla\BBB,\nabla\EEE)\,\frac{\partial\varphi}{\partial x_s}\,,
\nonumber
\end{align}
where $\fol_s$ is the polynomial expression above in $\BBB$, $\EEE$
and their derivatives.

Before proceeding we need the following technical lemma.

\begin{lemma}\label{perpp} If $X$ is a vector field on $M$ with values in
  $\R^{n+m}$, we have
\begin{align*}
(\nabla_i X)^\perp =
&\, \nabla_i X^\perp + \BBB_{ij}g^{js}\left\langle X,\partial\varphi/\partial
    x_s\right\rangle+\left\langle X,\BBB_{ij}\right\rangle
  g^{js}\frac{\partial\varphi}{\partial x_s}\,,\\
(\nabla_i X)^M=
&\, \nabla_i X^M - \BBB_{ij}g^{js}\left\langle X,\partial\varphi/\partial
    x_s\right\rangle-\left\langle X,\BBB_{ij}\right\rangle
  g^{js}\frac{\partial\varphi}{\partial x_s}\,.
\end{align*}
More in general,
\begin{align*}
(\nabla_{i_1\dots i_k} X)^\perp =
&\,\nabla_{i_1\dots i_k} X^\perp 
+\sum_{j=0}^{k-1}\pol^{k-j}_s(\nabla^{k-j-1}\BBB)\left\langle\nabla^jX,
\partial\varphi/\partial  x_s\right\rangle\\
&\,+\sum_{j=0}^{k-1}\left\langle \nabla^j
  X,\pol_s^{k-j}(\nabla^{k-j-1}\BBB)\right\rangle
\frac{\partial\varphi}{\partial x_s}\,.
\end{align*}
\end{lemma}
\begin{proof} We compute
\begin{align*}
(\nabla_iX)^\perp =
&\, \nabla_i X^\perp + (\nabla_i X^M)^\perp -
(\nabla_iX^\perp)^M\\
=&\,\nabla_i X^\perp + \BBB_{ij}g^{js}\left\langle X,\partial\varphi/\partial
    x_s\right\rangle - \left\langle\frac{\partial
    X^\perp}{\partial x_i},\frac{\partial\varphi}{\partial
    x_j}\right\rangle g^{js}\frac{\partial\varphi}{\partial x_s}\\
=&\,\nabla_i X^\perp + \BBB_{ij}g^{js}\left\langle X,\partial\varphi/\partial
    x_s\right\rangle+\left\langle X^\perp,\frac{\partial^2\varphi}{\partial
    x_i\partial x_j}\right\rangle g^{js}\frac{\partial\varphi}{\partial
  x_s}\\
=&\,\nabla_i X^\perp + \BBB_{ij}g^{js}\left\langle X,\partial\varphi/\partial
    x_s\right\rangle +\left\langle X,\BBB_{ij}\right\rangle
  g^{js}\frac{\partial\varphi}{\partial x_s}\,.
\end{align*}
The second formula is similar.\\
The third formula follows by induction, once one works in
normal coordinates where
$\frac{\partial^2\varphi}{\partial{x_i}\partial{x_j}}=\BBB_{ij}$,
which is a normal vector.
\end{proof}

\begin{rem}
Roughly, this lemma says that the interchange of differentiation and
projection operators introduces some extra terms in $\BBB$, $X$ and
their derivatives, and the order of differentiation of $X$ is
lower than the initial one.\\
This is useful when $X$ is a
function of $\BBB$, in particular, if $X$ is the mean curvature vector
$\HHH$ we have
\begin{align}\label{hperp}
(\nabla_{i_1\dots i_k} \HHH)^\perp =
&\,\nabla_{i_1\dots i_k} \HHH^\perp 
+\sum_{j=0}^{k-1}\pol^{k-j}_s(\nabla^{k-j-1}\BBB)\left\langle\nabla^j\HHH,
\partial\varphi/\partial  x_s\right\rangle\\
&\,+\sum_{j=0}^{k-1}\left\langle \nabla^j
  \HHH,\pol_s^{k-j}(\nabla^{k-j-1}\BBB)\right\rangle
\frac{\partial\varphi}{\partial x_s}\nonumber\\
=&\,\nabla_{i_1\dots i_k} \HHH
+\pol^k(\nabla^{k-1}\BBB) +\pol_s^k(\nabla^{k-1}\BBB)
\frac{\partial\varphi}{\partial x_s}\,,\nonumber
\end{align}
as $\HHH$ is a normal vector.
\end{rem}

\begin{lemma}\label{lemincasinato}
For any $s \in \mathbb{N}$ we have
\begin{align}\label{tort}
\frac{\partial\,}{\partial t}
\int_M\vert{\nabla^s\BBB}\vert^2\,d\mu
=&
\,-2\int_M \vert\nabla^{s+1}\BBB\vert^2\,d\mu
+\int_M\qol^{2s+4}(\nabla^{s+1}\BBB)\,d\mu
\\
&\,
-4\varepsilon k\int_M
\vert\nabla^{s+k-1}\BBB\vert^2\,d\mu
+\varepsilon\int_M\qol^{2k+2s}(\nabla^{2k+s-3}\BBB)\,d\mu\,,
\nonumber
\end{align}
and
\begin{equation}\label{rell}
\frac{\partial\,}{\partial t}
\int_M\vert{\BBB}\vert^{2s+2}\,d\mu
=\int_M\qol^{2s+4}(\nabla^{2}\BBB)\,d\mu + 
\varepsilon\int_M\qol^{2k+2s}(\nabla^{2k-3}\BBB)\,d\mu\,.
\end{equation}
\end{lemma}
\begin{proof}
Substituting ${\EEE}^\varepsilon$ in place of $\EEE$ in~\eqref{enzo} 
and expanding, after some computation using Lemma~\ref{perpp}, we get
\begin{align}
\derpar t {\BBB} 
= 
&\,2\varepsilon k(-1)^{k}
\nabla^{2}\Bigl(\overset{\text{$( k -2)$--times}}{\overbrace{\Delta^{M}\comp
\Delta^{M}\comp\,\dots\,\comp\Delta^{M}}}\,\,\HHH\Bigr)^\perp
+\nabla^2\HHH\label{ciao1000}
\\ 
&\,+\varepsilon\pol^{2k-1}(\nabla^{2k-3}\BBB)
+\varepsilon\pol_j^{2k-1}(\nabla^{2k-3}\BBB)\frac{\partial\varphi}{\partial
  x_j}+\pol_j^{3}(\nabla^{2}\BBB)\frac{\partial\varphi}{\partial
  x_j}\,.\nonumber
\end{align}
Applying now formula~\eqref{hperp} to the first term on the right
hand side of~\eqref{ciao1000}, we obtain
\begin{align*}
\derpar t {\BBB} = 
&\,2\varepsilon k(-1)^{k}
\nabla^{2}\,\,\overset{\text{$( k -2)$--times}}{\overbrace{\Delta^{M}\comp
\Delta^{M}\comp\,\dots\,\comp\Delta^{M}}}\,\,\HHH +\nabla^2\HHH\\ 
&\,+\varepsilon\pol^{2k-1}(\nabla^{2k-3}\BBB)
+\varepsilon\pol_j^{2k-1}(\nabla^{2k-3}\BBB)\frac{\partial\varphi}{\partial
  x_j}+\pol_j^{3}(\nabla^{2}\BBB)\frac{\partial\varphi}{\partial
  x_j}\,.
\end{align*} 

\noindent
We now observe that for any tensor $T$, we have
\begin{align*}
\frac{\partial\,}{\partial t} \nabla T = 
&\,\nabla
\frac{\partial\,}{\partial t} T 
+ T*\BBB*\nabla{\EEE}^\varepsilon
+T*\nabla\BBB*{\EEE}^\varepsilon\\
=&\, \nabla
\frac{\partial\,}{\partial t} T 
+ T*\pol^3(\nabla\BBB)
+\varepsilon T*\pol^{2k-1}(\nabla^{2k-4}\BBB)\,.
\end{align*}
Then, starting from equation~\eqref{ciao1000} and working by
induction, again in normal coordinates, we get
\begin{align*}
\frac{\partial\,}{\partial t} {\nabla^s\BBB} = 
&\,2\varepsilon k(-1)^{k}
\nabla^{s+2}\,\,\overset{\text{$( k -2)$--times}}{\overbrace{\Delta^{M}\comp
\Delta^{M}\comp\,\dots\,\comp\Delta^{M}}}\,\,\HHH + \nabla^{s+2}\HHH\\ 
&\,+\varepsilon\pol^{2k+s-1}(\nabla^{2k+s-3}\BBB)+\pol^{s+3}(\nabla^{s+1}\BBB)\\
&\,+\varepsilon\pol_j^{2k+s-1}(\nabla^{2k+s-3}\BBB)\frac{\partial\varphi}{\partial
  x_j}+\pol_j^{s+3}(\nabla^{s+2}\BBB)\frac{\partial\varphi}{\partial
  x_j}\,.
\end{align*}
Hence,
\begin{align*}
\frac{\partial\,}{\partial t} \vert{\nabla^s\BBB}\vert^2 = 
&\,4\varepsilon k(-1)^{k}
\nabla^{s+2}_{i_1\dots i_slw}\,\,\overset{\text{$( k -2)$--times}}{\overbrace{\Delta^{M}\comp
\Delta^{M}\comp\,\dots\,\comp\Delta^{M}}}\,\,\HHH\,
\nabla^s_{j_1\dots j_s}\BBB_{pz}\,g^{i_1j_1}\dots
g^{i_sj_s}g^{lp}g^{wz}\\ 
&\,+2\nabla^{s+2}_{i_1\dots i_slw}\HHH\,
\nabla^s_{j_1\dots j_s}\BBB_{pz}\,g^{i_1j_1}\dots
g^{i_sj_s}g^{lp}g^{wz}\\
&\,+\varepsilon\qol^{2k+2s}(\nabla^{2k+s-3}\BBB)+\qol^{2s+4}(\nabla^{s+1}\BBB)\,.
\end{align*}
Thus, we have that the time derivative of the quantity
$\int_M\vert{\nabla^s\BBB}\vert^2\,d\mu$ is given by 
\begin{align*}
&\,4\varepsilon k(-1)^{k}\int_M
\nabla^{s+2}_{i_1\dots i_slw}\,\,\overset{\text{$( k 
-2)$--times}}{\overbrace{\Delta^{M}\comp
\Delta^{M}\comp\,\dots\,\comp\Delta^{M}}}\,\,\HHH\,
\nabla^s_{j_1\dots j_s}\BBB_{pz}\,g^{i_1j_1}\dots 
g^{i_sj_s}g^{lp}g^{wz}\,d\mu
\\ 
&\,+ 2 \int_M \nabla^{s+2}_{i_1\dots i_slw}\HHH\,
\nabla^s_{j_1\dots j_s}\BBB_{pz}\,g^{i_1j_1}\dots 
g^{i_sj_s}g^{lp}g^{wz}\,d\mu
\\ 
&\,+\varepsilon\int_M\qol^{2k+2s}(\nabla^{2k+s-3}\BBB)\,d\mu
+\int_M\qol^{2s+4}(\nabla^{s+1}\BBB)\,d\mu\,,
\end{align*}
since 
$$
\frac{\partial\,}{\partial t}d\mu=\langle
  \HHH,{\EEE}^\varepsilon\rangle\,d\mu,
$$
hence its contribution
  can be absorbed in the last two terms.

We need now the following formula, which follows by direct computation,
\begin{equation}\label{maria}
\nabla_i\BBB_{jl}=\nabla_j\BBB_{il}+(\BBB*\BBB)_z
\frac{\partial\varphi}{\partial  x_z}\,.
\end{equation}
Indeed, the 3--tensor $\nabla_i\BBB_{jl}$ taking 
values in $\R^{n+m}$ is not symmetric (see also
Remark~\ref{rem3bbb}).\\
Reasoning then as in~\cite[Prop.~7.2-7.4]{mant5}, with the only care that
instead of applying Proposition 2.4 in that paper, we use 
formula~\eqref{maria},
we finally obtain~\eqref{tort} and~\eqref{rell}.
\end{proof}

\section{{$\varepsilon$}--Independent Estimates}\label{secepsindep}

For any integer $s > n/2$ and $\eps >0$ we set
\begin{equation}\label{defQepss}
Q^s_\eps(t)=\int_M 
\left(1+\vert\nabla^{s}\BBB\vert^2+\vert\BBB\vert^{2s+2}\right)\,d\mu, 
\qquad t \in [0,+\infty).
\end{equation}
Letting $\varepsilon>0$ vary, we want study the evolution of $Q_\eps^s$
under the flows associated with the 
functionals $\DG^\varepsilon$.\\

By the computations of the previous section we have
\begin{align}\label{derQeps}
\frac{\partial Q_\eps^s}{\partial t}
=&\,
-
2\int_M \vert\nabla^{s+1}\BBB\vert^2\,d\mu
-4\varepsilon k\int_M
\vert\nabla^{s+k-1}\BBB\vert^2\,d\mu
\\
&\,
+
\int_M\qol^{2s+4}(\nabla^{s+1}\BBB)\,d\mu
+
\varepsilon\int_M\qol^{2k+2s}(\nabla^{2k+s-3}\BBB)\,d\mu\,.
\nonumber
\end{align}

In order to deal with the polynomial terms we state in other words 
Proposition~6.5
in~\cite{mant5} (see all Section~6 in the same paper).

\begin{prop}\label{pnnp} Choosing some $\delta>0$ and 
setting $D={\mathrm {Vol}}(M)  +\Vert\HHH\Vert_{L^{n+\delta}(\mu)}$, 
there  exists a constant $C$ depending only on $n$, $m$, $l$, $z$, $j$, 
$p$,
  $q$, $r$, $\delta$ and $D$, such that for every manifold $(M,g)$ and
  covariant tensor $T=T_{i_1 \dots i_l}$, the following inequality holds
\begin{equation}\label{intercomp}
    \Vert\nabla^j T\Vert_{L^{p}{(\mu)}}\leq\,C\,\Vert 
T\Vert_{W^{z,q}{(\mu)}}^{\sigma}\Vert
      T\Vert_{L^r{(\mu)}}^{1-\sigma}\,,
\end{equation}
for all $z\in\NN$, $j\in\{0,\dots, z\}$, $p, q, r\in[1,+\infty)$ and $\sigma\in[j/z,1]$
with the compatibility condition
$$
\frac{1}{p}=\frac{j}{n}+\sigma\left(\frac{1}{q}-\frac{z}{n}\right)+\frac{1-\sigma}{r}\,.
$$
If such a condition gives a negative value for $p$, the inequality holds
in~\eqref{intercomp}
for every $p\in[1,+\infty)$ on the left hand side.
\end{prop}

This clearly implies, looking at the definition of the quantities
$Q_\eps^s$, that we can alternatively let the constant $C$ in 
inequality~\eqref{intercomp} depend on
$n$, $m$, $s$, $l$, $z$, $j$, $p$, $q$, $r$ and $Q_\eps^{[n/2]+1}$.

Working now as in Section~7 of~\cite{mant5}, with $s>n/2$ fixed, we can
interpolate the polynomial terms as follows,
\begin{align*}
\int_M\qol^{2s+4}(\nabla^{s+1}\BBB)\,d\mu\leq
&\,\int_M \vert\nabla^{s+1}\BBB\vert^2\,d\mu + C_1(Q_\eps^{[n/2]+1})
\\
\int_M\qol^{2k+2s}(\nabla^{2k+s-3}\BBB)\,d\mu\leq
&\,3k\int_M\vert\nabla^{s+k-1}\BBB\vert^2\,d\mu + C_2(Q_\eps^{[n/2]+1})
\end{align*}
where $C_1(Q_\eps^{[n/2]+1})$ and $C_2(Q_\eps^{[n/2]+1})$ are some 
constants depending only on 
$n$, $m$, $k$, $s$ and $Q_\eps^{[n/2]+1}$.
\\
Hence, for every $s>n/2$, by~\eqref{derQeps} we have the 
 estimate
\begin{equation*}
\frac{\partial Q_\eps^s}{\partial t}\leq
-\int_M \vert\nabla^{s+1}\BBB\vert^2\,d\mu
-\varepsilon k\int_M\vert\nabla^{s+k-1}\BBB\vert^2\,d\mu
+C_1(Q_\eps^{[n/2]+1})+\eps C_2(Q_\eps^{[n/2]+1})\leq 
C(Q_\eps^{[n/2]+1})\,,
\end{equation*}
where $C_1$, $C_2$ and $C$ depend on $\eps$ only through 
$Q_\eps^{[n/2]+1}$.  
\begin{teo}\label{teoteo3} For any integer $s>n/2$ there exists a smooth
  function $\ZZZ^s:\R\to(0,+\infty)$ such that
\begin{equation}\label{ciccio8}
\dert \int_M 
\left(1+\vert\nabla^{s}\BBB\vert^2+\vert\BBB\vert^{2s+2}\right)\,d\mu
\leq \ZZZ^s\left(\int_M
\left(  
1+\vert\nabla^{[n/2]+1}\BBB\vert^2+\vert\BBB\vert^{2[n/2]+4}\right)\,d\mu\right)
\end{equation}
for every $\eps \in (0,1)$ and any smooth evolution by the gradient of the
functional $\DG^\eps$.
\end{teo}
\begin{proof} The functions $\ZZZ^s$ can be clearly chosen to be
  smooth, possibly slightly enlarging the constants in the last inequality
  above.
\end{proof}

As a consequence we get the following corollary.
\begin{prop}\label{tos}
In the same hypotheses of Theorem~\ref{teoteo3}, 
there exists a continuous nonincreasing function $\Theta: (0, +\infty) \to 
(0, +\infty)$, 
{\em independent of $\varepsilon>0$}, such that for every $T\in\R$ and 
$t\in[T,T+\Theta(Q_\eps^{[n/2]+1}(T))]$ 
we have $Q_\eps^{[n/2]+1}(t)\leq
2Q_\eps^{[n/2]+1}(T)$.
\end{prop}
\begin{proof}
The statement follows by a standard ODE argument applied to the
differential inequality
\begin{equation*}
\dert \int_M 
\left(1+\vert\nabla^{[n/2]+1}\BBB\vert^2+\vert\BBB\vert^{2[n/2]+4}
\right)\,d\mu
\leq \ZZZ^{[n/2]+1}
\left(\int_M\left(  
1+\vert\nabla^{[n/2]+1}\BBB\vert^2+\vert\BBB\vert^{2[n/2]+4}\right)\,d\mu\right)
\end{equation*}
which is the first case of Theorem~\ref{teoteo3}.
\end{proof}

In other words, this proposition says that we have an
$\varepsilon$--independent control $Q_\eps^{[n/2]+1}\leq C$ in some
$\varepsilon$--independent time interval $[T,T+\Theta]$ 
(hence also a control the constants in Proposition~\ref{pnnp} 
and on the right hand side of inequalities~\eqref{ciccio8} for every
$s>n/2$), with $C$ and $\Theta$ depending (smoothly) only on
the value of $Q_\eps^{[n/2]+1}$ at the starting time $T$.

\section{Convergence to the Mean Curvature Flow}\label{seclim}

In this section we prove the convergence of solutions 
$\varphi^\eps:M\times[0,+\infty)\to\R^{n+m}$ to~\eqref{tttt1} 
(all starting from a common immersion $\varphi_0$) to the mean 
curvature flow $\varphi:M\times[0,T_{\mathrm {sing}})\to\R^{n+m}$ 
before its first singularity time. 

We need the following result (which can be proved as
in the codimension one case as in~\cite[Proposition~6.3]{mant5}).
\begin{prop}\label{pnnp2} If  a manifold $(M,g)$ satisfies
  ${\mathrm {Vol}}(M) +\Vert\HHH\Vert_{L^{n+\delta}(\mu)}\leq D$ for some
  $\delta>0$   then for every covariant tensor $S=S_{i_1 \dots i_l}$ we
  have
\begin{equation*}
\max_M\vert S\vert\leq C\left(\Vert\nabla S\Vert_{L^{p}(\mu)} 
+ \Vert S\Vert_{L^p(\mu)}\right)\qquad\text{ if $p>n$,}
\end{equation*}
where the constants depend only on $n$, $m$, $l$, $p$, $\delta$ and $D$.
\end{prop}

Let $Q^s(t)$ denote, for each nonnegative time $t$ before the first 
singularity, the right hand side of~\eqref{defQepss} for the mean 
curvature
flow $\varphi$ at time $t$.

\begin{lemma}\label{porol}
If the family of immersions  $\varphi^\eps(\cdot, T):M\to\R^{n+m}$ are 
bounded in the
$C^\infty$ topology, for any $s\in\NN$ 
all the quantities $\vert\nabla^{s}\BBB\vert$
are uniformly bounded by $\varepsilon$--independent
constants $C_s<+\infty$, in the time interval
$[T,T+\Theta(\sup_{\eps>0}Q^{[n/2]+1}_\eps(T))]$, where
$\Theta$ is the function in Proposition~\ref{tos}.
\end{lemma} 

\begin{proof}
By the $C^\infty$ boundedness of the family
$\varphi^\eps(\cdot, T):M\to\R^{n+m}$, all the quantities
$Q^{[n/2]+1}_\eps(T)$ are equibounded. 
As the function $\Theta$ is continuous and nonincreasing, setting
$\tau=\Theta(\sup_{\eps>0}Q^{[n/2]+1}_\eps(T))>0$, 
by Proposition~\ref{tos}
there exists a
constant $C>0$ such that $Q^{[n/2]+1}_\eps(t)\leq C$ for every
$\eps>0$ and $t\in[T,T+\tau]$.\\
Then, again by the boundedness of the family $\varphi^\eps(\cdot, T)$ and
Theorem~\ref{teoteo3}, in the same time interval $[T,T+\tau]$ all the  
quantities
\begin{equation*}
\int_M 
\left(1+\vert\nabla^{s}\BBB\vert^2+\vert\BBB\vert^{2s+2}\right)\,d\mu\,,
\end{equation*}
for every $s>n/2$, are bounded by $\varepsilon$--independent constants
$C_s<+\infty$. Moreover, all the constants in the interpolation
inequalities of Propositions~\ref{pnnp} and~\ref{pnnp2} are also bounded.\\
As a first step we see that, by means of
Proposition~\ref{pnnp}, we get the following estimates, for every
$p\in[2,+\infty)$ and $s\in\NN$, 
\begin{equation*}
\int_M \vert\nabla^{s}\BBB\vert^p\,d\mu\leq C_{s,p}
\end{equation*}
in the same time interval $[T,T+\tau]$. Here again the constants
$C_{s,p}<+\infty$ are $\varepsilon$--independent.\\
Then, we conclude by means of Proposition~\ref{pnnp2}.
\end{proof}

\begin{lemma}\label{piripi}
Assume that at time $t=T$ the family of maps
$\varphi^\varepsilon(\cdot,T):M\to\R^{n+m}$ converges as $\eps\to0$ in the
$C^\infty$ topology to the immersion $\varphi_T:M\to\R^{n+m}$.
Then the maps $\varphi^\eps$ smoothly converge in the time interval $[T, T+ 
\Theta(Q^{[n/2]+1}(T)))$ to the solution of the mean curvature flow starting from 
$\varphi_T$.
\end{lemma}

\begin{proof}
By the previous lemma, we have uniform bounds on $\BBB$ and its
derivatives in the time interval $[T,T+\tau]$ with 
$\tau=\Theta(\sup_{\eps>0}Q^{[n/2]+1}_\eps(T))$.
Then, there exists $C>0$ independent of $\eps$ such that 
$$
\left\vert\frac{\partial \varphi^\varepsilon(p,t)}{\partial t}\right\vert
=\vert\EEE^\eps(p,t)\vert< C
\qquad \forall (p,t)\in M\times[T,T+\tau], \  \eps>0.
$$
Now we consider the metric tensors 
$g_{ij}^\eps(p,t)=\left\langle\frac{\partial
    \varphi^\varepsilon(p,t)}{\partial x_i}, \frac{\partial
    \varphi^\varepsilon(p,t)}{\partial x_j}\right\rangle$, and 
fix a vector $V=\{v^i\}\in T_pM$. Then we have
$$
\left\vert 
\frac{\partial\,}{\partial t}\,\vert 
V\vert^2_{g^\eps(p,t)}\right\vert
=
\left\vert \partial_t
g_{ij}^\eps(p,t)v^iv^j\right\vert
=
2\left\vert \langle\EEE^\eps,\BBB_{ij}\rangle v^iv^j\right\vert
\leq
2\vert\EEE^\eps\vert\,\vert\BBB\vert_{g^\eps(p,t)}\vert
V\vert^2_{g^\eps(p,t)}
\leq 
C\vert V\vert^2_{g^\eps(p,t)}
$$
where $C$ does not depend on $\eps$.\\
Then a simple ODE argument shows  that the metrics $g^\eps_{ij}$ are
all equivalent; more precisely, there exists a positive constant $C$ 
depending only on $\varphi_T$ such that
\begin{equation}\label{metriceq}
\frac{{\mathrm {Id}}}{C}\leq g^\eps_{ij}(p,t)\leq C{\mathrm {Id}}\,,
\end{equation}
as matrices.\\
Moreover, as functions, all the 
$g^\eps_{ii}=\left\vert\frac{\partial\varphi^\eps}{\partial
    x_i}\right\vert^2$ are equibounded above by a common constant.\\
Hence, by Ascoli--Arzel\`a's Theorem, up to a subsequence, 
the immersions $\varphi^\eps$ uniformly converge, as $\eps\to0$ 
to some Lipschitz map $\widehat{\varphi}:M\times[T,T+\tau]\to\R^{n+m}$, which 
clearly satisfies $\widehat{\varphi}(p,T)=\varphi_T(p)$ for every 
$p\in M$.

Similarly, as the time derivative of the Christoffel symbols 
is given by
\begin{equation}\label{equ500}
\frac{\partial\,}{\partial t}\Gamma_{ij}^l = \nabla\EEE^\eps * {\BBB} +
\EEE^\eps * \nabla\BBB
\end{equation}
(see the beginning of Section~\ref{evolgeosec}) and all the metrics
are equivalent, it follows that all the Christoffel symbols are
equibounded. This means that 
estimating the covariant derivatives is equivalent
to estimate the standard derivatives in coordinates, hence, we have
immediately $\vert\partial^s\nabla^l\BBB\vert\leq C_{s,l}$ 
for every $s,l\in\NN$.\\
Since
$$
\frac{\partial\,}{\partial t}g_{ij}^\eps = 2\langle\EEE^\eps,\BBB_{ij}\rangle
$$
we get
$$
\left\vert\nabla^s\frac{\partial\,}{\partial
t}g_{ij}\right\vert \leq C_s\,,
$$
and, by formula~\eqref{equ500}, 
$$
\left\vert\nabla^s\frac{\partial\,}{\partial
t}\Gamma_{ij}^l\right\vert\leq C_s\,.
$$
for every $s\in\NN$.\\
Hence, we get
$\left\vert\partial^s\frac{\partial\,}{\partial
  t}\Gamma_{ij}^l\right\vert\leq C_s$  which implies, as the family of
maps $\varphi^\eps_T$ is bounded in the $C^\infty$--topology, that 
$\vert\partial^s\Gamma_{ij}^l\vert\leq C_s$.\\
Since we already know that $\vert\varphi^\eps\vert$ are equibounded, 
$\vert \partial\varphi^\eps\vert \leq C$ and 
${\partial^2\varphi^\eps}=\Gamma\partial\varphi^\eps+\BBB$, by the 
estimates $\vert\partial^s\nabla^l\BBB\vert\leq C_{s,l}$, 
we can conclude that the derivatives $\vert\partial^s\varphi^\eps\vert$ 
are all bounded by $\eps$--independent constants $C_s$, for every
$s\in\NN$.\\
Finally, the uniform control on the time and mixed derivatives of
$\varphi^\eps$ follows using the evolution equation.

Hence, the sub--convergence $\varphi^\eps\to\widehat{\varphi}$, 
as $\eps\to 0$, is in the $C^{\infty}$ topology and $\widehat \varphi$ is
smooth, moreover, the limit metric is positive definite
by~\eqref{metriceq}.\\
Passing to the limit in the evolution equation
$\partial_t\varphi^\eps=\EEE^\eps$, by the bounds on $\BBB$ and its
derivatives, shows that
$\widehat{\varphi}:M\times[T,T+\tau]\to\R^{n+m}$ 
flows by mean curvature with a starting smooth datum $\varphi_T$. Since
this flow is unique, {\em all} the sequence of maps $\varphi^\eps$ converges to
$\widehat{\varphi}$ which hence coincides with $\varphi$.

Chosen now any $\delta>0$, let $\eps_0>0$ be such that
$$
\sup_{0<\eps<\eps_0}Q^{[n/2]+1}_\eps(T)-Q^{[n/2]+1}(T)<\delta\,.
$$
Since $\Theta$ is nonincreasing (see
Lemma~\ref{porol}), in the interval
$[T,T+\Theta(Q^{[n/2]+1}(T)+\delta)]$ the sequence $\varphi^\eps$
converges to the mean curvature flow $\varphi$.
Letting $\delta$ to zero, as the function $\Theta$ is continuous, we get
the thesis.
\end{proof}

We are now in the position to conclude the proof of  Theorem~\ref{mainth}. 

\begin{proof}[Proof of Theorem~\ref{mainth}]
Let $T_{\mathrm {max}}$ be the maximal time
such that $\varphi^\eps$ converge to the solution of the mean 
curvature flow equation $\varphi$ in
$C^\infty(M\times [0,T_{\mathrm {max}}))$ starting at time $0$ from the common
immersion $\varphi_0$. Observe that $T_{\mathrm {max}}$ is positive 
by Lemma~\ref{piripi}. We want to show that $T_{\mathrm {max}}$ coincides
with the first singularity time $T_{\mathrm {sing}}$ for $\varphi$.

Assume by contradiction that $T_{\mathrm {max}} < T_{\mathrm {sing}}$.
Then $\varphi(\cdot,t) \to \varphi(\cdot, {T_{\mathrm {max}}})$ in 
$C^\infty(M)$ as $t\to T_{\mathrm
  {max}}$.
As the function $\Theta$ is continuous, there exists
$$
\lim_{t\to T_{\mathrm  {max}}}
\Theta(Q^{[n/2]+1}(t))=\Theta(Q^{[n/2]+1}({T_{\mathrm
    {max}}}))=\theta>0\,.
$$
Choosing now a time $T\in[T_{\mathrm  {max}}-\theta/4,T_{\mathrm
  {max}})$ such that $\Theta(Q^{[n/2]+1}(T))>\theta/2$, and applying
Lemma~\ref{piripi}, we see that $\varphi^\eps(\cdot, t)$ converges to the 
mean curvature flow also for $t$ in the interval $[T, T+\theta/2]$. 
As $T+\theta/2>T_{\mathrm  {max}}-\theta/4+\theta/2>T_{\mathrm
  {max}}$, we have a contradiction.
\end{proof}

\bibliographystyle{amsplain}
\bibliography{biblio}

\end{document}